\pdfoutput=1

\documentclass[10pt]{article}

\usepackage{url}
\usepackage{mathtools}
\usepackage{amssymb}
\usepackage{amsthm}
\usepackage{latexsym}
\usepackage[shortlabels]{enumitem}
\usepackage{dsfont}
\usepackage{appendix}
\usepackage{color} 
\usepackage[utf8]{inputenc}
\usepackage[T1]{fontenc}
\usepackage{geometry}
\usepackage{todonotes}
\usepackage{lmodern}
\usepackage{anyfontsize}
\usepackage{stmaryrd}
\usepackage{comment}
\usepackage{bm}
\usepackage{mathrsfs}
\usepackage{mdframed}
\usepackage{natbib}
\usepackage[english]{babel}
\usepackage{cases}

\mdfsetup{middlelinecolor=blue,middlelinewidth=2pt,linewidth=0pt,backgroundcolor=blue!15,,roundcorner=10pt}

\allowdisplaybreaks

\setcitestyle{numbers,open={[},close={]}}

\definecolor{red}{rgb}{0.7,0.15,0.15}
\definecolor{green}{rgb}{0,0.5,0}
\definecolor{blue}{rgb}{0,0,0.7}

\usepackage{hyperref}
\hypersetup{colorlinks, linkcolor={red}, citecolor={green}, urlcolor={blue}}
\numberwithin{equation}{section}

\newtheorem{theorem}{Theorem}[section]
\newtheorem{assumption}[theorem]{Assumption}

\newtheorem{lemma}[theorem]{Lemma}

\newtheorem{proposition}[theorem]{Proposition}

\newtheorem{remark}[theorem]{Remark}

\usepackage[capitalize,nameinlink,noabbrev]{cleveref}

\crefname{assumption}{assumption}{assumptions}
\Crefname{assumption}{Assumption}{Assumptions}
\crefname{assumptio}{assumption}{assumptions}
\Crefname{assumptio}{Assumption}{Assumptions}
\crefname{assumptionalt}{assumption}{assumptions}
\Crefname{assumptionalt}{Assumption}{Assumptions}
\crefname{lemma}{lemma}{lemmas}
\Crefname{lemma}{Lemma}{Lemmas}
\crefname{claim}{claim}{claims}
\Crefname{claim}{Claims}{Claims}
\crefname{proposition}{proposition}{propositions}
\Crefname{proposition}{Proposition}{Propositions}
\crefname{definition}{definition}{definitions}
\Crefname{definition}{Definition}{Definitions}
\crefname{remark}{remark}{remarks}
\Crefname{remark}{Remark}{Remarks}

\AddToHook{env/assumption/begin}{\crefalias{theorem}{assumption}}
\AddToHook{env/lemma/begin}{\crefalias{theorem}{lemma}}
\AddToHook{env/proposition/begin}{\crefalias{theorem}{proposition}}
\AddToHook{env/definition/begin}{\crefalias{theorem}{definition}}
\AddToHook{env/remark/begin}{\crefalias{theorem}{remark}}
\AddToHook{env/claim/begin}{\crefalias{theorem}{claim}}

\newcommand\bu{\mathbf{u}}
\newcommand\bv{\mathbf{v}}
\newcommand\bw{\mathbf{w}}

\newcommand\cA{\mathcal A}

\newcommand\cE{\mathcal E}
\newcommand\cF{\mathcal F}

\newcommand\cI{\mathcal I}

\def \E{\mathbb{E}}
\def \F{\mathbb{F}}

\def \L{\mathbb{L}}

\def \N{\mathbb{N}}
\def \P{\mathbb{P}}

\def \R{\mathbb{R}}

\def \W{\mathbb{W}}

\def\bnothing{{\color{blue}\pmb{\varnothing}}}

\def\bone{{\color{blue}\mathbf{1}}}
\def\btwo{{\color{blue}\mathbf{2}}}

\def\bd{{\color{blue}\mathbf{d}}}

\def\bi{{\color{blue}\mathbf{i}}}
\def\bj{{\color{blue}\mathbf{j}}}

\def\bm{{\color{blue}\mathbf{m}}}
\def\bn{{\color{blue}\mathbf{n}}}

\def\bu{{\color{blue}\mathbf{u}}}
\def\bv{{\color{blue}\mathbf{v}}}
\def\bw{{\color{blue}\mathbf{w}}}

\newcommand{\1}{\mathbf{1}}
\newcommand{\shuprod}{\mathrel{\sqcup \mkern -3.2mu \sqcup}}

\def\d{\mathrm{d}}

\title{Signature representations of certain BSDEs
\footnotetext{
Support from the Deutsche Forschungsgemeinschaft	CRC/TRR 388	Project ID 516748464 is gratefully acknowledged.}
}
\date{November 12, 2026}

\author{Peter {\sc Bank}\footnote{TU Berlin, Germany, bank@math.tu-berlin.de.}
\and Christian {\sc Bayer}\footnote{Weierstrass Institute, Germany, christian.bayer@wias-berlin.de.}
\and Marco {\sc Rauscher}\footnote{TU Munich, Germany, marco.rauscher@tum.de.}
\and Marco {\sc Rodrigues}\footnote{TU Berlin and Weierstrass Institute, Germany, marco.rodrigues@tu-berlin.de and marco.rodrigues@wias-berlin.de.}
\and Nikolas {\sc Tapia}\footnote{HU Berlin and Weierstrass Institute, Germany, tapia@wias-berlin.de.}
}

\begin{document}

\maketitle

\begin{abstract}
We link the solvability of backward stochastic differential equations with time-augmented Brownian signature data to the solvability of a differential equation on the tensor algebra of signature coefficients. We show that in certain cases the solutions are of closed form and remain within the class of linear functionals of the signature. Lastly, we derive a closed-form solution to a pricing and hedging problem in terms of the time-augmented Brownian signature.
\end{abstract}

\section{Introduction}

\medskip
Backward stochastic differential equations (BSDEs) were systematically introduced by \citeauthor*{pardoux1990adapted} \cite{pardoux1990adapted,pardoux1992backward}. Ever since their introduction, they have been applied to stochastic control problems, partial differential equations, principal-agent problems, mathematical finance, nonlinear expectation, numerical methods and machine learning; see \citeauthor*{el1997backward} \cite{el1997backward} for a survey of the relation of BSDEs to stochastic control and mathematical finance and \citeauthor*{zhang2017backward} \cite{zhang2017backward} for a survey of the general methodology, subsequent developments, and further references.
One of their main advantages is that they naturally allow for the modelling of path-dependent and non-Markovian dynamics. Existence and uniqueness of solutions to BSDEs are often obtained by fixed-point arguments or monotone approximation schemes. However, having applications in mind, one is also interested in obtaining closed-form solutions, often a challenging task depending on the exact data of the BSDE under consideration.

\medskip
One way of encoding path dependence in stochastic equations is through the path signature---originating in the work of \citeauthor*{chen1957integration} \cite{chen1957integration,chen1977iterated} and subsequently generalized by \citeauthor*{lyons1998differential} \cite{lyons1998differential}---which is the sequence of all iterated path integrals of the driving noise.
Since its introduction, the signature has been used in the development of rough path theory, applied to numerical methods for SDEs, to machine learning, to stochastic control, and to mathematical finance, to name a few; we refer to the recent edited volume by \citeauthor*{bayer2026signature} \cite{bayer2026signature}.

\medskip
Our main motivation comes from the application of the signature methodology in the work of \citeauthor*{abi2024path} \cite{abi2024path}. They are able to represent certain non-Markovian stochastic differential equations as linear (differential) equations on the tensor algebra of signature coefficients. Subsequently solving the linear equation yields a solution to the original SDE as a linear functional of the signature. Here, we adopt a similar methodology for BSDEs. More precisely, consider the following path-dependent BSDE
\begin{equation}\label{eq_intro_bsde}
    Y_t = \xi + \int_t^T f_s(Y_s,Z_s)\d s - \int_t^T Z_s\d W_s, \; t \in [0,T], \; \textnormal{$\P$--a.s.},
\end{equation}
where $W$ denotes a Brownian motion. A natural question arises: under what assumptions on the terminal condition $\xi$ and the generator $f$ does the solution $(Y,Z)$ have a representation in terms of the time-augmented $\widehat{W}_t = (t,W_t)$ Brownian signature $\widehat{\W}$?

\medskip
Suppose that $Y_t = \langle \boldsymbol{\eta}(t),\widehat{\W}_t\rangle$ is a linear functional obtained by pairing the tensor $\boldsymbol{\eta}(t)$ with the signature $\widehat{\W}_t$. Formally, we obtain the dynamics
\[
    \bigg\langle \dot{\boldsymbol{\eta}}(t) + \boldsymbol{\eta}(t)|_\bone + \frac{1}{2}\boldsymbol{\eta}(t)|_{\btwo\btwo}, \widehat{\W}_t(\omega) \bigg\rangle + f\Big(t,\omega,\langle \boldsymbol{\eta}(t),\widehat{\W}_t\rangle, \big\langle \boldsymbol{\eta}(t)|_\btwo, \widehat{\W}_t(\omega) \big\rangle\Big) = 0, \; \langle \boldsymbol{\eta}(T),\widehat{\W}_T(\omega)\rangle = \xi(\omega).
\]
If the generator $f$ and terminal condition $\xi$ satisfy the following algebraic conditions
\begin{equation}\label{eq_polynomial_generator}
    f\Big(t,\omega,\langle \boldsymbol{\eta}(t),\widehat{\W}_t\rangle, \big\langle \boldsymbol{\eta}(t)|_\btwo, \widehat{\W}_t(\omega) \big\rangle\Big) 
    = \big\langle \mathbf{f}^{\boldsymbol{\eta}(\cdot)}(t),\widehat{\W}_t(\omega) \big\rangle \; \textnormal{and} \; \xi(\omega) = \langle \boldsymbol{\zeta}, \widehat{\W}_T(\omega)\rangle,
\end{equation}
for some signature coefficients $\mathbf{f}^{\boldsymbol{\eta}(\cdot)}(t)$ and $\boldsymbol{\zeta}$, it seems appropriate to consider the following tensor differential equation for the signature coefficients $\big\{\boldsymbol{\eta}(t) \,|\, t \in [0,T]\big\}$:
\begin{equation}\label{eq_tensor_differential_equation}
\begin{cases}
    \displaystyle \dot{\boldsymbol{\eta}}(t) + \boldsymbol{\eta}(t)\big|_{\bone} + \frac{1}{2}\boldsymbol{\eta}(t)\big|_{\btwo\btwo} + \mathbf{f}^{\boldsymbol{\eta}(\cdot)}(t)
    = 0, \\[0.5em]
    \boldsymbol{\eta}(T) = \boldsymbol{\zeta}.
\end{cases}
\end{equation}
There are two main objectives in this work. First, we derive a verification result that relates solutions to the tensor differential equation \eqref{eq_tensor_differential_equation} to solutions of the BSDE \eqref{eq_intro_bsde}. Secondly, we compute solutions to \eqref{eq_tensor_differential_equation} in closed form for an affine generator. When the generator is identically zero, this corresponds to the martingale representation problem. Here, we can transform the time-dependent tensor $\boldsymbol{\eta}(t)$ into a time-independent tensor $\bar{\boldsymbol{\eta}}$, which is interesting in and of itself.

\medskip
Related differential equations for tensors have already appeared in previous work, for example, in \citeauthor*{abi2025signature} \cite{abi2025signature}, \citeauthor*{jaber2026affine} \cite{jaber2026affine}, and \citeauthor*{cuchiero2026signature} \cite{cuchiero2026signature} to compute conditional expectations and transforms, and in \citeauthor*{jaber2025signature} \cite{jaber2025signature}, and \citeauthor*{abi2026stochastic} \cite{abi2026stochastic} in relation to stochastic control and hedging.

\medskip
The remainder is organised as follows. In \Cref{sec_preliminaries}, we describe the setting and introduce signatures. \Cref{sec_verification} contains the main verification result, which connects the tensor differential equation to BSDEs. In \Cref{sec_examples}, we solve the tensor differential equation in the case of BSDEs with affine generators and obtain a closed-form solution of the BSDE. Moreover, we also obtain a martingale representation result with a time-independent tensor equation. \Cref{sec_application} contains an application to pricing and hedging. Lastly, in \Cref{sec_further}, we provide some concluding remarks.

\section{Preliminaries and setting}\label{sec_preliminaries}

Let $(\Omega,\cF,\P)$ be a probability space supporting a (one-dimensional) Brownian motion $W=(W_t)_{t \in [0,\infty)}$, that is, $W$ is a real-valued, right-continuous and $\P$--a.s. continuous process starting at zero with independent and stationary increments such that $\E^\P[W^2_t] = t$ for every $t \in [0,\infty)$. Let $\F = (\cF_t)_{t \in [0,\infty)}$ be the filtration generated by $W$, and let $T \in [0,\infty)$ denote a finite time horizon.

\begin{remark}
    The restriction to one-dimensional Brownian motion is solely for notational simplicity$;$ the multidimensional case is analogous. 
\end{remark}

\subsection{BSDEs}\label{sec_pre_bsdes}

Throughout this work, we fix a real-valued and $\cF_T$-measurable function $\xi$ and a map 
\[
    f : [0,T]\times\Omega\times\R\times\R \longrightarrow \R, \; (t,\omega,y,z) \longmapsto f(t,\omega,y,z),
\]
jointly measurable with respect to the progressive $\sigma$-algebra on $[0,T]\times\Omega$ and the Borel $\sigma$-algebra on $\R\times\R$.
A solution to the BSDE with data $(f,\xi)$ is any pair of processes $(Y,Z)$ consisting of two real-valued, $\F$-progressive processes $Y = (Y_t)_{t \in [0,T]}$ and $Z = (Z_t)_{t \in [0,T]}$ such that 
\begin{equation}\label{eq_bsde_lipschitz}
\begin{cases}
    \displaystyle \int_0^T |f_s(Y_s,Z_s)|\d s + \int_0^T |Z_s|^2 \d s < \infty, \; \textnormal{$\P$--a.s.}, & \\[1em]
    \displaystyle Y_t = \xi + \int_t^T f_s(Y_s,Z_s)\d s - \int_t^T Z_s \d W_s, \; t \in [0,T], \; \textnormal{$\P$--a.s.} &
\end{cases}
\end{equation}
It is well known (see \citeauthor*{zhang2017backward} \cite[Section~4.3                                                                                                      ]{zhang2017backward}) that in case $\xi$ is $\P$--square-integrable,
$f(t,\omega,\cdot,\cdot)$ is Lipschitz, uniformly in $(t,\omega)$, and
\[
    \E^\P\Bigg[\bigg(\int_0^T |f_s(0,0)|\d s\bigg)^2\Bigg] < \infty,
\]
there exists a unique square-integrable solution $(Y,Z)$ to the BSDE \eqref{eq_bsde_lipschitz} with data $(f,\xi)$, that is,
\begin{equation}\label{eq_square_integrability}
    \E^\P\bigg[\sup_{t \in [0,T]}|Y_t|^2 + \int_0^T |Z_s|^2 \d s \bigg] < \infty.
\end{equation}
The pair $(Y,Z)$ is unique in the following sense: any other pair $(Y^\prime,Z^\prime)$ of real-valued, $\F$-progressive processes $Y^\prime =(Y^\prime_t)_{t \in [0,T]}$ and $Z^\prime=(Z^\prime_t)_{t \in [0,T]}$ such that \eqref{eq_bsde_lipschitz}--\eqref{eq_square_integrability} hold satisfies $Y=Y^\prime$ up to $\P$-indistinguishability and $Z=Z^\prime$ outside a $\P\otimes\d s$--null set.

\subsection{Signature notations}

We fix some tensor-algebraic notation, which we borrow from \cite{jaber2025signature}. Here, $\N$ denotes the non-negative integers and $\N^\ast \coloneqq \N\setminus\{0\}$.

\subsubsection{Tensor algebras and words}

Fix $d \in \N^\ast$. For $n \in \N$, write $(\R^d)^{\otimes n}$ for the $n$-fold tensor power of $\R^d$; we use the convention $(\R^d)^{\otimes 0}\coloneqq\R$. The extended tensor algebra over $\R^d$ is
\[
    T((\R^d))
    \coloneqq
    \{ \boldsymbol{\ell}=(\boldsymbol{\ell}^n)_{n\in\N} \,|\, \boldsymbol{\ell}^n\in (\R^d)^{\otimes n} \text{ for all } n\in\N \}.
\]
The tensor algebra $T(\R^d)$ is defined as
\[
    T(\R^d) \coloneqq \bigcup_{M\in\N} T^M(\R^d), \; \textnormal{where} \; T^M(\R^d) \coloneqq \{ \boldsymbol{\ell}\in T((\R^d)) \,|\, \boldsymbol{\ell}^n=0 \text{ for all } n>M \}.
\]
It is convenient to write $\boldsymbol{\ell} \in T((\R^d))$ formally as a sum as follows. Let $\{e_1,\dots,e_d\}$ be the canonical basis of $\R^d$. We identify $e_{i_1}\otimes\cdots\otimes e_{i_n}$ with the word $\mathbf{\bi_\bone\cdots \bi_\bn}$ over the alphabet $\{\bone,\dots,\bd\}$, and the unit $1\in\R$ with the empty word $\bnothing$. Set
\[
    V\coloneqq \bigcup_{n\in\N}V_n,
\]
where $V_n$ denotes the collection of words of length $n$; $V_0$ only contains the empty word $\bnothing$. Any $\boldsymbol{\ell}\in T((\R^d))$ is formally written as
\[
    \boldsymbol{\ell} = \sum_{\bv\in V} \boldsymbol{\ell}^{\bv}\bv,
\]
where $\boldsymbol{\ell}^{\bv}$ denotes the coefficient of $\boldsymbol{\ell}$ corresponding to the word $\bv$. In what follows, we use the notation $\bv$ and $\bw$ for words, and denote the length of $\bv$ by $|\bv|$.

\medskip
For $\boldsymbol{\ell},\boldsymbol{p} \in T((\R^d))$, the tensor product $\boldsymbol{\ell}\otimes \boldsymbol{p}$ is defined by
\[
    (\boldsymbol{\ell}\otimes \boldsymbol{p})^\bw \coloneqq \sum_{\substack{\bu,\bv \in V \\ \bu\bv=\bw}} \boldsymbol{\ell}^\bu \boldsymbol{p}^\bv, \, \bw \in V.
\]
We use the conventions
\[
    \boldsymbol{\ell}^{\otimes 0} \coloneqq \bnothing \; \textnormal{and} \; \boldsymbol{\ell}^{\otimes n} \coloneqq \underbrace{\boldsymbol{\ell}\otimes\cdots\otimes\boldsymbol{\ell}}_{\textnormal{$n$-times}} \; \textnormal{for $n \in \N^\ast$.}
\]

\subsubsection{The tensor-product exponential}

For $\boldsymbol{\ell}\in T((\R^d))$, let $\exp_\otimes(\boldsymbol{\ell})\in T((\R^d))$ be given by
\[
    \exp_\otimes(\boldsymbol{\ell})^\bw \coloneqq \sum_{k=0}^\infty\frac{(\boldsymbol{\ell}^{\otimes k})^\bw}{k!}
\]
This is well-defined for the following reason. Write $\boldsymbol{\ell} = \boldsymbol{\ell}^\bnothing\bnothing + \tilde{\boldsymbol{\ell}}$. Then
\[
    (\boldsymbol{\ell}^{\otimes k})^\bw = \sum_{j = 0}^{\min\{k,|\bw|\}} \binom{k}{j}(\boldsymbol{\ell}^\bnothing)^{k-j}(\tilde{\boldsymbol{\ell}}^{\otimes j})^\bw
\]
implies
\[
    \sum_{k=0}^\infty\frac{|(\boldsymbol{\ell}^{\otimes k})^\bw|}{k!} < \infty.
\]
Moreover,
\[
    \exp_\otimes(\boldsymbol{\ell}) = \mathrm{e}^{\boldsymbol{\ell}^\bnothing} \exp_\otimes(\tilde{\boldsymbol{\ell}}).
\]

\subsubsection{The right-shift operators \texorpdfstring{$\boldsymbol{\ell}\vert_{\bw}$}{l restricted to w}
and \texorpdfstring{$\boldsymbol{\ell}\vert_{\boldsymbol{\zeta}}$}{l restricted to zeta}}

For $\boldsymbol{\ell}\in T((\R^d))$ and $\bw\in V$, we define $\boldsymbol{\ell}|_{\bw} \in T((\R^d))$ by
\[
    (\boldsymbol{\ell}|_\bw)^\bv \coloneqq \boldsymbol{\ell}^{\bv\bw}, \; \bv \in V,
\]
and for $\boldsymbol{\zeta} \in T((\R^d))$ with $\sum_{\bw \in V} |\boldsymbol{\ell}^{\bv\bw}\boldsymbol{\zeta}^\bw| < \infty$ for all $\bv \in V$, we define $\boldsymbol{\ell}|_{\boldsymbol{\zeta}}\in T((\R^d))$ by
\[
    (\boldsymbol{\ell}|_{\boldsymbol{\zeta}})^\bv \coloneqq  \sum_{\bw\in V} \boldsymbol{\ell}^{\bv\bw}\boldsymbol{\zeta}^{\bw}, \; \bv \in V.
\]
In particular, $\boldsymbol{\ell}|_{\boldsymbol{\zeta}}$ exists if $\boldsymbol{\ell} \in T(\R^d)$ or $\boldsymbol{\zeta} \in T(\R^d)$.

\subsubsection{Shuffle products}

A shuffle of $\bv=\bv_\bone\cdots\bv_\bm$ and $\bw=\bw_\bone\cdots\bw_\bn$, where each $\bv_\bi$ and $\bw_\bj$ is a letter in the alphabet $\{1,\ldots,d\}$, is any word $\bu$ of length $m+n$ obtained by interlacing the letters of $\bv$ and $\bw$ in such a way that the order of the $\bv_\bi$'s and the order of the $\bw_\bj$'s are both preserved. The shuffle product $\bv \shuprod \bw \in T((\R^d))$ of $\bv,\bw \in V$ is defined by
\[
    \bv \shuprod \bw \coloneqq \sum_{\bu \in \operatorname{Sh}(\bv,\bw)}N(\bv,\bw;\bu) \bu,
\]
where $\operatorname{Sh}(\bv,\bw)$ denotes the set of all distinct shuffles of $\bv$ and $\bw$, and $N(\bv,\bw;\bu)$ counts how many interlacings of letters of $\bv$ and $\bw$ produce the word $\bu$. 
Note that $\bnothing \shuprod \bv = \bv \shuprod \bnothing = \bv$. For $\boldsymbol{\ell},\boldsymbol{\zeta} \in T((\R^d))$, we define $\boldsymbol{\ell}\shuprod\boldsymbol{\zeta} \in T((\R^d))$ by
\[
    (\boldsymbol{\ell}\shuprod \boldsymbol{\zeta})^{\bu} = \sum_{\substack{\bv,\bw\in V\\ \bu\in \operatorname{Sh}(\bv,\bw)}} N(\bv,\bw;\bu) \boldsymbol{\ell}^{\bv}\boldsymbol{\zeta}^{\bw}, \; \bu \in V.
\]

\subsubsection{Time-augmented Brownian signature}

Let $\widehat{W} = (\widehat{W}_t)_{t \in [0,\infty)}$ be the $\R^2$-valued, time-augmented Brownian motion $\widehat{W}_t \coloneqq (t,W_t)$. The Stratonovich signature of $\widehat{W}$ is the map
\[
    \widehat{\W}: \Omega \times [0,\infty) \ni (\omega,t) \longmapsto (1,\widehat{\W}^1_t(\omega),\ldots,\widehat{\W}^n_t(\omega),\ldots)\in T((\R^2)),
\]
defined recursively by $\widehat{\W}^0_t(\omega)\coloneqq 1$ and
\[
    \widehat{\W}^{m+1}_t(\omega) \coloneqq \bigg(\int_0^t \widehat{\W}^m_s\otimes \circ \d \widehat{W}_s\bigg)(\omega), \; m\in\N,
\]
where the stochastic integral is taken componentwise in the Stratonovich sense. For $\boldsymbol{\ell} \in T((\R^2))$ and $t \in [0,T]$, we then define
\[
    \langle \boldsymbol{\ell},\widehat{\mathbb W}_t\rangle = \sum_{n = 0}^\infty\sum_{\substack{\bw\in V \\ |\bw| = n}} \boldsymbol{\ell}^\bw \widehat{\W}_t^\bw,
\]
provided the (potentially infinite) sum is well-defined; this is the case if, for example,
\[
    \|\boldsymbol{\ell}\|^{\cA(\widehat{\W})}_t \coloneqq \sum_{n = 0}^\infty \bigg| \sum_{\substack{\bw \in V \\ |\bw|=n}} \boldsymbol{\ell}^\bw \widehat{\W}_t^\bw\bigg| < \infty, \; \textnormal{$\P$--a.s.}
\]
In particular, if $\|\boldsymbol{\ell}\|^{\cA(\widehat{\W})}_t < \infty$, then $|\langle \boldsymbol{\ell}, \widehat{\W}_t \rangle|\leq \|\boldsymbol{\ell}\|^{\cA(\widehat{\W})}_t$.

\medskip
For $\boldsymbol{\ell},\boldsymbol{\zeta} \in T((\R^2))$ and $t \in [0,T]$ with $\|\boldsymbol{\ell}\|^{\cA(\widehat{\W})}_t + \|\boldsymbol{\zeta}\|^{\cA(\widehat{\W})}_t < \infty$, $\P$--a.s., the pairing $\langle \boldsymbol{\ell}\shuprod \boldsymbol{\zeta},\widehat{\mathbb W}_t\rangle$ is $\P$--a.s. well-defined, and the shuffle product yields
\[
    \langle \boldsymbol{\ell},\widehat{\mathbb W}_t\rangle \langle \boldsymbol{\zeta},\widehat{\mathbb W}_t\rangle
    = \langle \boldsymbol{\ell}\shuprod \boldsymbol{\zeta},\widehat{\mathbb W}_t\rangle, \; \textnormal{$\P$--a.s.}
\]

\subsubsection{It\^o formula}

Let $\cI^\prime(\widehat{\W})$ denote the collection of all $\boldsymbol{\eta} : [0,T] \longrightarrow T((\R^2))$ such that $\|\boldsymbol{\eta}(t)\|^{\cA(\widehat{\W})}_t < \infty$, $\P$--a.s. for every $t \in [0,T]$, each $\boldsymbol{\eta}^\bw$ is absolutely continuous and
\[
    \int_0^T \Big(\|\dot{\boldsymbol{\eta}}(s)\|^{\cA(\widehat{\W})}_s + \|\boldsymbol{\eta}(s)|_{\bone}\|^{\cA(\widehat{\W})}_s + \|\boldsymbol{\eta}(s)|_{\btwo\btwo}\|^{\cA(\widehat{\W})}_s + \big(\|\boldsymbol{\eta}(s)|_{\btwo}\|^{\cA(\widehat{\W})}_s\big)^2\Big) \d s < \infty, \; \textnormal{$\P$--a.s.},
\]
where $\dot{\boldsymbol{\eta}} = \sum_{\bw\in V} \dot{\boldsymbol{\eta}}^\bw \bw$ and $\dot{\boldsymbol{\eta}}^\bw$ denotes the density of $\boldsymbol{\eta}^\bw$ with respect to Lebesgue measure.

\begin{lemma}\label{lem_ito_signature}
    For any $\boldsymbol{\eta} \in \cI^\prime(\widehat{\W})$, the pairings 
    \[
        \langle \boldsymbol{\eta}(s)|_{\bone}, \widehat{\W}_s\rangle, \; \langle \boldsymbol{\eta}(s)|_{\btwo\btwo},\widehat{\W}_s\rangle, \; \langle  \dot{\boldsymbol{\eta}}(s),\widehat{\W}_s\rangle, \; \textnormal{and} \; \langle \boldsymbol{\eta}(s)|_{\btwo},\widehat{\W}_s\rangle,
    \]
    are \textnormal{$\P\otimes\d s$--a.e.} well-defined and admit $\P\otimes\d s$-representatives that are $\F$-progressive and satisfy
    \[
        \int_0^T \big( |\langle  \dot{\boldsymbol{\eta}}(s),\widehat{\W}_s\rangle| + |\langle \boldsymbol{\eta}(s)|_{\bone}, \widehat{\W}_s\rangle| + |\langle \boldsymbol{\eta}(s)|_{\btwo\btwo},\widehat{\W}_s\rangle| + |\langle \boldsymbol{\eta}(s)|_{\btwo},\widehat{\W}_s\rangle|^2 \big) \d s < \infty, \; \textnormal{$\P$--a.s.}
    \]
    Moreover,
    \begin{equation}\label{eq_ito_signatures}
        \langle \boldsymbol{\eta}(t), \widehat{\W}_t\rangle 
        = \boldsymbol{\eta}(0)^{\bnothing} + \int_0^t \Big\langle \dot{\boldsymbol{\eta}}(s) + \boldsymbol{\eta}(s)|_{\bone} + \frac{1}{2} \boldsymbol{\eta}(s)|_{\btwo\btwo},\widehat{\W}_s\Big\rangle \d s + \int_0^t \langle \boldsymbol{\eta}(s)|_{\btwo},\widehat{\W}_s\rangle\d W_s, \; \textnormal{$\P$--a.s.}, \; t \in [0,T].
    \end{equation}
\end{lemma}

\begin{remark}
    The formula \eqref{eq_ito_signatures} appeared in \textnormal{\cite[Theorem 3.3 and Corollary 3.4]{abi2024path}} and we include a proof here because an integrability condition on $\|\dot{\boldsymbol{\eta}}(s)\|^{\cA(\widehat{\W})}_s$ was not explicitly stated. Without such a condition, the term $\langle \dot{\boldsymbol{\eta}}(s),\widehat{\W}_s \rangle$ may not be integrable with respect to Lebesgue measure.
\end{remark}

\begin{proof}
    For $M \in \N$, let
    \[
        \boldsymbol{\eta}^M(t) \coloneqq \sum_{n = 0}^M\sum_{\substack{\bw \in V \\ |\bw| = n}} \boldsymbol{\eta}^\bw(t) \bw, \; t \in [0,T].
    \]
    Then
    \[
        \langle\boldsymbol{\eta}^M(t),\widehat{\W}_t\rangle = \sum_{n = 0}^M \sum_{\substack{\bw \in V \\ |\bw| = n}} \boldsymbol{\eta}^\bw(t)\widehat{\W}_t^\bw.
    \]
    Since each $\boldsymbol{\eta}^\bw$ is absolutely continuous, an application of the integration by parts formula yields
    \[
        \boldsymbol{\eta}^\bw(t)\widehat{\W}_t^\bw = \boldsymbol{\eta}^\bw(0)\widehat{\W}_0^{\bw} + \int_0^t \boldsymbol{\eta}^\bw(s)\d\widehat{\W}_s^{\bw} + \int_0^t \dot{\boldsymbol{\eta}}^\bw(s) \widehat{\W}_s^{\bw} \d s, \; t \in [0,T], \; \textnormal{$\P$--a.s.}
    \]
    By definition of the signature, we have
    \[
        \d\widehat{\W}_s^{\bv\bone} = \widehat{\W}_s^{\bv} \d s
        \; \textnormal{and} \;
        \d\widehat{\W}_s^{\bv\btwo} = \widehat{\W}_s^{\bv} \circ \d W_s = \widehat{\W}_s^{\bv} \d W_s + \frac{1}{2}\d[\widehat{\W}^{\bv},W]_s , \; \textnormal{$\P$--a.s.}
    \]
    If $\bv = \bnothing$ or $\bv = \bu \bone$, then $\d [\widehat{\W}^{\bv},W] \equiv 0$, and if $\bv = \bu \btwo$, then $\d [\widehat{\W}^{\bv},W] = \d[\widehat{\W}^{\bu \btwo},W] = \widehat{\W}^{\bu}_s\d s$. Hence
    \[
        \langle\boldsymbol{\eta}^M(t),\widehat{\W}_t\rangle
        = \boldsymbol{\eta}^\bnothing(0) +  \int_0^t\Big\langle \boldsymbol{\eta}^M(s)|_\bone + \frac{1}{2} \boldsymbol{\eta}^M(s)|_{\btwo\btwo} + \dot{\boldsymbol{\eta}}^M(s), \widehat{\W}_s\Big\rangle\d s  + \int_0^t \big\langle \boldsymbol{\eta}^M(s)|_\btwo,\widehat{\W}_s\rangle \d W_s, \; t \in [0,T], \; \textnormal{$\P$--a.s.}
    \]
    Define
    \[
    \langle \boldsymbol{\eta}(s)|_\bone, \widehat{\W}_s\rangle \coloneqq H_s \1_{\{H_s \in \R\}}, \; \textnormal{where} \; H_s \coloneqq \limsup_{M\rightarrow\infty} \langle \boldsymbol{\eta}^M(s)|_\bone, \widehat{\W}_s\rangle, \; s \in [0,T].
    \]
    Then, $\langle \boldsymbol{\eta}(\cdot)|_\bone, \widehat{\W}_\cdot\rangle$ is $\F$-progressive, and since $\int_0^T \|\boldsymbol{\eta}(s)|_\bone\|^{\cA(\widehat{\W})}_s \d s < \infty$, $\P$--a.s., we have
    \[
        \langle \boldsymbol{\eta}(s)|_\bone, \widehat{\W}_s\rangle = \lim_{M \rightarrow\infty} \langle \boldsymbol{\eta}^M(s)|_\bone, \widehat{\W}_s\rangle, \; \textnormal{$\d s$--a.e. and in $\L^1(\d s)$, $\P$--a.s.}
    \]
    We define the terms $\langle \boldsymbol{\eta}(s)|_{\btwo\btwo},\widehat{\W}_s\rangle,$ $\langle  \dot{\boldsymbol{\eta}}(s),\widehat{\W}_s\rangle,$ and $\langle \boldsymbol{\eta}(s)|_{\btwo},\widehat{\W}_s\rangle$ analogously. The truncations $\langle \boldsymbol{\eta}^M(s)|_{\btwo\btwo},\widehat{\W}_s\rangle$ and $\langle  \dot{\boldsymbol{\eta}}^M(s),\widehat{\W}_s\rangle$ then also converge to $\langle \boldsymbol{\eta}(s)|_{\btwo\btwo},\widehat{\W}_s\rangle$ and $\langle  \dot{\boldsymbol{\eta}}(s),\widehat{\W}_s\rangle$ $\d s$--a.e. and in $\L^1(\d s)$, $\P$--a.s., respectively. The truncation $\big\langle \boldsymbol{\eta}^M(s)|_\btwo,\widehat{\W}_s\rangle$ instead converges to $\big\langle \boldsymbol{\eta}(s)|_\btwo,\widehat{\W}_s\rangle$ in $\L^2(\d s)$, $\P$--almost surely. 
    The claim \eqref{eq_ito_signatures} now follows from an application of dominated convergence for Lebesgue integrals and from continuity of the stochastic integral (see \cite[Proposition~18.6]{kallenberg2021foundations}). This completes the proof.
\end{proof}

\begin{remark}\label{rem_progressive_representative}
The preceding result also provides an $\F$-progressive, right-continuous and \textnormal{$\P$--a.s.} continuous, $\P$-modifi\-cation of  $\langle \boldsymbol{\eta}(\cdot), \widehat{\W}_\cdot\rangle$, which we use from now on.
\end{remark}

\section{Main verification result}\label{sec_verification}

The main verification result is based on a structural property of the generator, which we formulate as an assumption.

\begin{assumption}\label{ass_algebraic_property}
For any $\boldsymbol{a}\in\cI^\prime(\widehat{\W})$, there exists $\mathbf{f}^{\boldsymbol{a}(\cdot)} : [0,T] \longrightarrow T((\R^{2}))$ such that 
\begin{enumerate}
    \item[$(i)$] $\|\mathbf{f}^{\boldsymbol{a}(\cdot)}(s)\|^{\cA(\widehat{\W})}_{s}<\infty$, \textnormal{$\P\otimes\d s$--a.e.}, and
    \item[$(ii)$]\label{eq_structural_condition_generator} $f_s\big(\langle \boldsymbol{a}(s),\widehat{\W}_s\rangle,\langle \boldsymbol{a}(s)|_{\btwo}, \widehat{\W}_s\rangle\big) 
    = \big\langle \mathbf{f}^{\boldsymbol{a}(\cdot)}(s),\widehat{\W}_s \big\rangle$, $\P\otimes \d s$--a.e.
\end{enumerate}
\end{assumption}

This assumption allows us to relate solutions to a deterministic differential equation on the tensor space $T((\R^2))$ to solutions to the BSDE \eqref{eq_bsde_lipschitz}. The following verification result makes this precise.

\begin{theorem}\label{thm_main_verification}
Suppose that the generator $f$ satisfies \textnormal{\Cref{ass_algebraic_property}}, and that $\xi = \langle \boldsymbol{\zeta} , \widehat{\W}_T\rangle$ for $\boldsymbol{\zeta} \in T((\R^{2}))$ satisfying $\|\boldsymbol{\zeta}\|^{\cA(\widehat{\W})}_T < \infty$, \textnormal{$\P$--almost surely}.
If there exists $\boldsymbol{\eta} \in \cI^\prime(\widehat{\W})$ satisfying
\begin{equation}\label{eq_coefficient_equation}
\begin{cases}
\displaystyle \dot{\boldsymbol{\eta}}(s) + \boldsymbol{\eta}(s)\big|_{\bone} + \frac{1}{2}\boldsymbol{\eta}(s)\big|_{\btwo\btwo} + \mathbf{f}^{\boldsymbol{\eta}(\cdot)}(s)
= 0,
& \textnormal{$\d s$--a.e. on $[0,T]$}, \\[0.5em]
\boldsymbol{\eta}(T) = \boldsymbol{\zeta},
\end{cases}
\end{equation}
then $(Y,Z)$ given by
\[
    Y_t \coloneqq \langle \boldsymbol{\eta}(t),\widehat{\W}_t\rangle \; \textnormal{and} \; Z_t \coloneqq \langle \boldsymbol{\eta}(t)|_{\btwo},\widehat{\W}_t\rangle, \; \textnormal{$t \in [0,T]$},
\]
is a solution to the \textnormal{BSDE \eqref{eq_bsde_lipschitz}}.
\end{theorem}

\begin{remark}
    \begin{enumerate}
        \item[$(i)$] The condition $\boldsymbol{\eta} \in \cI^\prime(\widehat{\W})$ ensures that all relevant integrals are well-defined. 
        \item[$(ii)$] This verification result does not assert existence nor uniqueness of a solution to the tensor differential equation \eqref{eq_coefficient_equation}, and the differential equation may, in general, be truly nonlinear. 
        Establishing existence of a solution hinges on the dependence of $\mathbf{f}^{\boldsymbol{\eta}(\cdot)}$ on $\boldsymbol{\eta}(\cdot)$, and we expect, as in classical \textnormal{BSDE} theory $($see \textnormal{\citeauthor*{bao2011backward} \cite{bao2011backward}}, and \textnormal{\citeauthor*{xing2012on} \cite[Equation~0.2]{xing2012on}}$)$, that uniqueness depends on the class in which a solution is sought, in general.
    \end{enumerate}
\end{remark}

\begin{proof}
By \Cref{lem_ito_signature}, the assumption on the generator $f$, and the coefficient equation for $\boldsymbol{\eta}$,
\begin{align*}
    \langle \boldsymbol{\eta}(t), \widehat{\W}_t\rangle 
    &= \boldsymbol{\eta}(0)^{\bnothing} + \int_0^t \Big\langle \boldsymbol{\eta}(s)|_{\bone} + \frac{1}{2} \boldsymbol{\eta}(s)|_{\btwo\btwo} + \dot{\boldsymbol{\eta}}(s),\widehat{\W}_s\Big\rangle \d s + \int_0^t \langle \boldsymbol{\eta}(s)|_{\btwo},\widehat{\W}_s\rangle\d W_s \\
    &= \boldsymbol{\eta}(0)^{\bnothing} - \int_0^t \big\langle \mathbf{f}^{\boldsymbol{\eta}(\cdot)}(s),\widehat{\W}_s \big\rangle \d s + \int_0^t \langle \boldsymbol{\eta}(s)|_{\btwo},\widehat{\W}_s\rangle\d W_s\\
    &= \boldsymbol{\eta}(0)^{\bnothing} - \int_0^t f_s\big(\langle \boldsymbol{\eta}(s),\widehat{\W}_s\rangle,\langle \boldsymbol{\eta}(s)|_{\btwo}, \widehat{\W}_s\rangle\big) \d s + \int_0^t \langle \boldsymbol{\eta}(s)|_{\btwo},\widehat{\W}_s\rangle\d W_s, \; \textnormal{$t\in[0,T]$, $\P$--a.s.}
\end{align*}
It follows from \Cref{rem_progressive_representative} that $Y$ admits a right-continuous and $\F$-progressive $\P$-modification and from \Cref{lem_ito_signature} that $Z$ admits an $\F$-progressive $\P\otimes\d s$-representative. We use these versions below. Moreover, $\boldsymbol{\eta} \in \cI^\prime(\widehat{\W})$, \Cref{ass_algebraic_property}.$(ii)$, and \eqref{eq_coefficient_equation} imply that
\[
    \int_0^T \big(|Z_s|^2  + |f_s(Y_s,Z_s)|\big) \d s < \infty, \; \textnormal{$\P$--a.s.}
\]
Since $Y_T = \langle \boldsymbol{\eta}(T),\widehat{\W}_T\rangle = \langle \boldsymbol{\zeta}, \widehat{\W}_T\rangle = \xi$, $\P$--a.s., it follows that
\[
    Y_t = \xi + \int_t^T f_s(Y_s,Z_s) \d s - \int_t^T Z_s \d W_s, \; t \in [0,T], \; \textnormal{$\P$--a.s.}
\]
Hence $(Y,Z)$ is a solution to the BSDE \eqref{eq_bsde_lipschitz}. This completes the proof.
\end{proof}

\section{Examples}\label{sec_examples}

We illustrate the verification result of \Cref{thm_main_verification} in two concrete settings.

\subsection{BSDEs with affine generator}\label{sec_affine_bsdes}

Consider a generator of the form 
\[
    f(s,\omega,y,z) = \alpha(s)y + \beta(s)z + \langle\boldsymbol{\gamma}(s),\widehat{\W}_s(\omega)\rangle, \; (s,\omega,y,z) \in [0,T]\times\Omega\times\R\times\R.
\]
The tensor differential equation \eqref{eq_coefficient_equation} turns into the following inhomogeneous linear differential equation:
\begin{equation}\label{eq_coefficient_equation_deterministic}
\begin{cases}
    \dot{\boldsymbol{\eta}}(s) + \boldsymbol{\eta}(s)\big|_{\bone} + \frac{1}{2}\boldsymbol{\eta}(s)\big|_{\btwo\btwo} + \alpha(s) \boldsymbol{\eta}(s) +  \beta(s) \boldsymbol{\eta}(s)\big|_{\btwo} + \boldsymbol{\gamma}(s)
    = 0,
    & \text{$\d s$--a.e. on $[0,T]$}, \\[0.5em]
    \boldsymbol{\eta}(T) = \boldsymbol{\zeta},
\end{cases}
\end{equation}
where the role of $\mathbf{f}^{\boldsymbol{\eta}(\cdot)}$ is played by
\[
    \mathbf{f}^{\boldsymbol{\eta}(\cdot)}(s) = \alpha(s) \boldsymbol{\eta}(s) + \beta(s) \boldsymbol{\eta}(s)\big|_{\btwo} + \boldsymbol{\gamma}(s).
\]
The candidate solution in the following result was obtained from a variation of constants ansatz; the particular homogeneous case $\alpha = \beta = \boldsymbol{\gamma} = 0$ was considered in \cite{jaber2025signature,jaber2026affine}.

\begin{proposition}\label{prop_affine_generator}
Let $\boldsymbol{\zeta} \in T(\R^2)$, let $\alpha,\beta : [0,T] \longrightarrow \R$ be Borel-measurable with $\alpha$ Lebesgue-integrable and $\beta$ bounded, and let $\boldsymbol{\gamma} : [0,T] \longrightarrow T^M(\R^2)$, for some $M \in \N^\ast$, be Borel-measurable and Lebesgue-integrable componentwise. Suppose that  
\[
    \xi = \langle \boldsymbol{\zeta}, \widehat{\W}_T\rangle \; \textnormal{and} \; f(s,\omega,y,z) = \alpha(s)y + \beta(s)z + \langle \boldsymbol{\gamma}(s),\widehat{\W}_s\rangle, \; (s,\omega,y,z) \in [0,T]\times\Omega\times\R\times\R.
\]
Let
\begin{align*}
    \cE (t,s) &\coloneqq  \exp\bigg(\int_t^s \alpha(r) \d r\bigg) \bigg(\bnothing + \sum_{n = 1}^\infty \int_{t < t_1 <\cdots < t_n < s} g(t_1)\otimes \cdots \otimes g(t_n) \d t_1 \ldots \d t_n\bigg), \; \textnormal{where} \; g(r) \coloneqq \bone + \frac{1}{2}\btwo\btwo + \beta(r)\btwo.
\end{align*}
Then
\begin{equation}\label{eq_eta_alpha_beta}
    \boldsymbol{\eta}(t) \coloneqq \int_t^T \boldsymbol{\gamma}(s)|_{\cE(t,s)} \d s + \boldsymbol{\zeta}|_{\cE(t,T)}, \; t \in [0,T],
\end{equation}
is in $\cI^\prime(\widehat{\W})$ and is a solution to \eqref{eq_coefficient_equation_deterministic}. Consequently, $(Y,Z)$ given by
\[
    Y_t = \langle \boldsymbol{\eta}(t),\widehat{\W}_t \rangle, \; Z_t = \langle \boldsymbol{\eta}(t)|_\btwo, \widehat{\W}_t\rangle, \; t \in [0,T], \; \textnormal{$\P$--a.s.},
\]
is a solution to the \textnormal{BSDE \eqref{eq_bsde_lipschitz}}. If, additionally, $\alpha$ is bounded, then $(Y,Z)$ is the unique square-integrable solution.
\end{proposition}

\begin{proof}
Since $\boldsymbol{\zeta} \in T(\R^2)$, we suppose without loss of generality that $\boldsymbol{\zeta},\boldsymbol{\gamma}(t) \in T^{M}(\R^2)$ for all $t \in [0,T]$; otherwise we increase $M$. This then implies that $\boldsymbol{\eta}(t) \in T^M(\R^2)$ for all $t \in [0,T]$. Clearly, $\boldsymbol{\eta}(T) = \boldsymbol{\zeta}$. We note that
\begin{equation}\label{eq_differential_E}
    \partial_t\cE(t,s) = - \bigg( \bone + \frac{1}{2}\btwo\btwo + \alpha(t)\bnothing + \beta(t)\btwo \bigg)\otimes \cE(t,s), \; \textnormal{$\d t$--a.e. on $[0,s]$, for every $s \in [0,T]$}.
\end{equation}
By definition, there exists $C \in (0,\infty)$ such that, for all $\bw \in V$ with $|\bw| \leq M$,
\begin{align}\label{eq_fubini_integrability}
    &|\cE(t,s)^\bw| \leq C, \; 0 \leq t \leq s \leq T, \; \textnormal{and} \nonumber\\
    &|\partial_t\cE(t,s)^\bw| \leq C (1+|\alpha(t)| + |\beta(t)|), \; \textnormal{$\d t$--a.e. on $[0,s]$, for every $s \in [0,T]$}.
\end{align}
This implies also that, for any $\bu,\bw \in V$ with $|\bw|\leq M$,
\[
    \int_0^T\int_0^s |\boldsymbol{\gamma}(s)^{\bu}\partial_t\cE(t,s)^\bw| \d t \, \d s < \infty.
\]
Hence, the components in the finite sum
\begin{align}\label{eq_pairing_bold_eta_linear_ode}
    \boldsymbol{\eta}(t)^\bv
    = \sum_{\substack{\bw\in V}}\bigg(\int_t^T \boldsymbol{\gamma}(s)^{\bv \bw} \cE(t,s)^\bw \d s + \boldsymbol{\zeta}^{\bv\bw}\cE(t,T)^\bw \bigg)
\end{align}
satisfy, by Fubini's theorem,
\begin{align*}
    \int_t^T \boldsymbol{\gamma}(s)^{\bv \bw} \cE(t,s)^\bw \d s 
    &= \int_t^T \boldsymbol{\gamma}(s)^{\bv \bw} \bigg( \cE(s,s)^\bw - \int_t^s \partial_u \cE(u,s)^\bw \d u\bigg) \d s \\
    &= \int_t^T \boldsymbol{\gamma}(s)^{\bv \bw} \cE(s,s)^\bw \d s - \int_t^T \boldsymbol{\gamma}(s)^{\bv \bw}\int_t^s \partial_u \cE(u,s)^\bw \d u \d s \\
    &= \int_t^T \boldsymbol{\gamma}(s)^{\bv \bw} \cE(s,s)^\bw \d s - \int_t^T \int_u^T \boldsymbol{\gamma}(s)^{\bv \bw} \partial_u \cE(u,s)^\bw \d s \d u .
\end{align*}
This immediately yields absolute continuity of the components in $\boldsymbol{\eta}(t)$. With \eqref{eq_differential_E}, we obtain
\begin{align*}
    \dot{\boldsymbol{\eta}}(t)^\bv 
    &= \sum_{\bw\in V} \bigg( - \boldsymbol{\gamma}(t)^{\bv\bw} \cE(t,t)^{\bw} + \int_t^T \boldsymbol{\gamma}(s)^{\bv\bw}\partial_t\cE(t,s)^{\bw}\d s + \boldsymbol{\zeta}^{\bv\bw}\partial_t\cE(t,T)^\bw\bigg) \\
    &= -\boldsymbol{\gamma}(t)^\bv -  \sum_{\bw\in V} \int_t^T \boldsymbol{\gamma}(s)^{\bv\bw}\Bigg(\bigg( \bone + \frac{1}{2}\btwo\btwo + \alpha(t)\bnothing + \beta(t)\btwo \bigg)\otimes \cE(t,s) \Bigg)^\bw\d s \\
    &\quad -\sum_{\bw\in V}\boldsymbol{\zeta}^{\bv\bw} \Bigg(\bigg( \bone + \frac{1}{2}\btwo\btwo + \alpha(t)\bnothing + \beta(t)\btwo \bigg)\otimes \cE(t,T)\Bigg)^\bw \\
    &= -\boldsymbol{\gamma}(t)^\bv - \sum_{\bw\in V}\int_t^T \boldsymbol{\gamma}(s)^{\bv\bw} \Bigg( \sum_{\substack{\bu\in V \\ \bone\bu=\bw}}\cE(t,s)^\bu + \frac{1}{2}\sum_{\substack{\bu\in V \\ \btwo\btwo\bu=\bw}}\cE(t,s)^\bu + \alpha(t)\cE(t,s)^\bw + \beta(t)\sum_{\substack{\bu\in V \\ \btwo\bu=\bw}}\cE(t,s)^\bu \Bigg)\d s \\
    &\quad - \sum_{\bw \in V} \boldsymbol{\zeta}^{\bv\bw} \Bigg( \sum_{\substack{\bu \in V \\ \bone\bu=\bw}}\cE(t,T)^\bu+ \frac{1}{2}\sum_{\substack{\bu\in V \\ \btwo\btwo\bu=\bw}}\cE(t,T)^\bu + \alpha(t)\cE(t,T)^\bw + \beta(t)\sum_{\substack{\bu\in V \\ \btwo\bu=\bw}}\cE(t,T)^\bu \Bigg) \\
    &= -\boldsymbol{\gamma}(t)^\bv - \sum_{\bw\in V}\int_t^T \Bigg( \sum_{\substack{\bu\in V \\ \bone\bu=\bw}} \boldsymbol{\gamma}(s)^{\bv\bone\bu}\cE(t,s)^\bu \\
    &\quad + \frac{1}{2}\sum_{\substack{\bu\in V \\ \btwo\btwo\bu=\bw}} \boldsymbol{\gamma}(s)^{\bv\btwo\btwo\bu}\cE(t,s)^\bu  + \alpha(t)\boldsymbol{\gamma}(s)^{\bv\bw}\cE(t,s)^\bw + \beta(t)\sum_{\substack{\bu\in V \\ \btwo\bu=\bw}} \boldsymbol{\gamma}(s)^{\bv\btwo\bu}\cE(t,s)^\bu \Bigg)\d s \\
    &\quad - \sum_{\bw\in V} \Bigg( \sum_{\substack{\bu \in V \\ \bone\bu=\bw}}\boldsymbol{\zeta}^{\bv\bone\bu}\cE(t,T)^\bu + \frac{1}{2}\sum_{\substack{\bu \in V \\ \btwo\btwo\bu=\bw}}\boldsymbol{\zeta}^{\bv\btwo\btwo\bu}\cE(t,T)^\bu + \alpha(t)\boldsymbol{\zeta}^{\bv\bw}\cE(t,T)^\bw + \beta(t)\sum_{\substack{\bu\in V \\ \btwo\bu=\bw}}\boldsymbol{\zeta}^{\bv\btwo\bu}\cE(t,T)^\bu \Bigg) \\
    &= -\boldsymbol{\gamma}(t)^\bv - \sum_{\bw\in V}\int_t^T \Bigg( \boldsymbol{\gamma}(s)^{\bv\bone\bw}\cE(t,s)^\bw + \frac{1}{2}\boldsymbol{\gamma}(s)^{\bv\btwo\btwo\bw}\cE(t,s)^\bw + \alpha(t)\boldsymbol{\gamma}(s)^{\bv\bw}\cE(t,s)^\bw + \beta(t)\boldsymbol{\gamma}(s)^{\bv\btwo\bw}\cE(t,s)^\bw \Bigg)\d s \\
    &\quad - \sum_{\bw\in V} \Bigg( \boldsymbol{\zeta}^{\bv\bone\bw}\cE(t,T)^\bw + \frac{1}{2}\boldsymbol{\zeta}^{\bv\btwo\btwo\bw}\cE(t,T)^\bw + \alpha(t)\boldsymbol{\zeta}^{\bv\bw}\cE(t,T)^\bw + \beta(t)\boldsymbol{\zeta}^{\bv\btwo\bw}\cE(t,T)^\bw \Bigg) \\
    &= -\boldsymbol{\gamma}(t)^\bv  - \sum_{\bw\in V}\int_t^T \Bigg(\boldsymbol{\gamma}(s)^{\bv\bone\bw}+ \frac{1}{2}\boldsymbol{\gamma}(s)^{\bv\btwo\btwo\bw}+ \alpha(t)\boldsymbol{\gamma}(s)^{\bv\bw} + \beta(t)\boldsymbol{\gamma}(s)^{\bv\btwo\bw} \Bigg)\cE(t,s)^\bw \d s \\
    &\quad - \sum_{\bw\in V} \Bigg( \boldsymbol{\zeta}^{\bv\bone\bw} + \frac{1}{2}\boldsymbol{\zeta}^{\bv\btwo\btwo\bw} + \alpha(t)\boldsymbol{\zeta}^{\bv\bw} + \beta(t)\boldsymbol{\zeta}^{\bv\btwo\bw} \Bigg)\cE(t,T)^\bw \\
    &= - \bigg(\boldsymbol{\gamma}(t)^\bv + \boldsymbol{\eta}^{\bv\bone}(t) + \frac{1}{2}\boldsymbol{\eta}^{\bv\btwo\btwo}(t) + \alpha(t)\boldsymbol{\eta}^\bv(t) + \beta(t)\boldsymbol{\eta}^{\bv\btwo}(t) \bigg), \; \textnormal{$\d t$--a.e. on $[0,T]$}.
\end{align*}
Since $\boldsymbol{\gamma}$ and $\alpha$ are Lebesgue-integrable, $\beta$ is bounded, and $\boldsymbol{\eta}$ is componentwise absolutely continuous, we have
\[
    \int_0^T \|\dot{\boldsymbol{\eta}}(s)\|^{\cA(\widehat{\W})}_s \d s < \infty, \; \textnormal{$\P$--a.s.}
\]
Since \eqref{eq_pairing_bold_eta_linear_ode} is a finite sum, $\|\boldsymbol{\eta}(t)\|^{\cA(\widehat{\W})}_t < \infty$ for all $t \in [0,T]$, and since $\boldsymbol{\eta}(t)$ is componentwise continuous, it follows that 
\[
    \int_0^T \Big( \|\boldsymbol{\eta}(s)|_{\bone}\|^{\cA(\widehat{\W})}_s + \|\boldsymbol{\eta}(s)|_{\btwo\btwo}\|^{\cA(\widehat{\W})}_s + \big(\|\boldsymbol{\eta}(s)|_{\btwo}\|^{\cA(\widehat{\W})}_s\big)^2\Big) \d s < \infty, \; \textnormal{$\P$--a.s.}
\]
Hence, $\boldsymbol{\eta}\in \cI^\prime(\widehat{\W})$ and it solves the differential equation \eqref{eq_coefficient_equation_deterministic}.
That $(Y,Z)$ is a solution to the BSDE \eqref{eq_bsde_lipschitz} follows from \Cref{thm_main_verification}. This proves the first assertion.

\medskip
We turn to the second assertion and now additionally assume that $\alpha$ is bounded; this immediately implies that the generator $f$ is uniformly Lipschitz continuous.
Note that the components of $\widehat{\W}$ have square-integrable running suprema on $[0,T]$. Thus, $\xi$ is $\P$--square-integrable, as $\boldsymbol{\zeta}$ is a finite tensor, and from
\[
    \E^\P\Bigg[\bigg(\int_0^T |\langle \boldsymbol{\gamma}(s),\widehat{\W}_s\rangle| \d s\bigg)^2\Bigg]^{1/2} \leq \sum_{\substack{\bw \in V \\ |\bw| \leq M}} \E^\P\bigg[\sup_{s \in [0,T]}|\widehat{\W}_s^{\bw}|^2\bigg]^{1/2} \int_0^T |\boldsymbol{\gamma}(s)^\bw|\d s < \infty,
\]
we obtain
\[
    \E^\P\Bigg[\bigg(\int_0^T |f_s(0,0)|\d s\bigg)^2\Bigg] < \infty.
\]
Hence, there exists a unique square-integrable solution to the BSDE \eqref{eq_bsde_lipschitz}; see \cite[Theorem~4.3.1]{zhang2017backward}. 
Since $\boldsymbol{\eta}(\cdot)$ is $T^M(\R^2)$-valued and componentwise absolutely continuous, it follows that $(Y,Z)$ satisfies \eqref{eq_square_integrability} and is therefore the unique square-integrable solution. This completes the proof.
\end{proof}

\begin{remark}
    Under suitable integrability assumptions, the $Y$-component of a linear \textnormal{BSDE} has a representation as a conditional expectation under a probability measure equivalent to $\P$$;$ see \textnormal{\cite[Section~4.1]{zhang2017backward}}. By contrast, the previous result provides an explicit representation for both $Y$ and $Z$ in terms of the time-augmented Brownian signature.
\end{remark}

\subsection{Martingale representation}

We consider the special affine generator $f \equiv 0$ and a terminal condition of the form $\xi = \langle \boldsymbol{\zeta}, \widehat{\W}_T\rangle$ for $\boldsymbol{\zeta} \in T(\R^2)$. Solving the corresponding BSDE amounts to finding the appropriate Brownian integrand in the martingale representation of
\[
    Y_t \coloneqq \E^\P\big[\langle \boldsymbol{\zeta}, \widehat{\W}_T\rangle\big|\cF_t\big], \; t \in [0,T].
\]
We thus refer to \eqref{eq_bsde_lipschitz} as the martingale representation problem in this section.
The relevant tensor differential equation becomes
\begin{equation}
\begin{cases}
    \dot{\boldsymbol{\eta}}(s) + \boldsymbol{\eta}(s)\big|_{\bone} + \frac{1}{2}\boldsymbol{\eta}(s)\big|_{\btwo\btwo} = 0,
    & \text{$\d s$--a.e. on $[0,T]$}, \\[0.5em]
    \boldsymbol{\eta}(T) = \boldsymbol{\zeta},
\end{cases}
\end{equation}
and by \Cref{prop_affine_generator}, its solution is
\[
    \boldsymbol{\eta}(t)=\boldsymbol{\zeta}|_{\cE(t,T)}, \; \cE(t,T) = \exp_\otimes\bigg((T-t)\Big(\bone+\frac{1}{2}\btwo\btwo\Big)\bigg).
\]
Nonetheless, it turns out that $Y_t = \langle \boldsymbol{\eta}(t),\widehat{\W}_t\rangle = \langle \bar{\boldsymbol{\eta}}, \widehat{\W}_t\rangle$, \textnormal{$\P\otimes\d t$--a.e.}, for some $\bar{\boldsymbol{\eta}} \in T(\R^2)$ that can be constructed from $\{\boldsymbol{\eta}(t) \,|\, t \in [0,T]\}$. The following gives a precise statement; compare with the representations in \cite[Theorem~4.2]{jaber2025signature} and \cite[Proposition~5.8]{jaber2026malliavin}.

\begin{proposition}
Let $\boldsymbol{\zeta} \in T(\R^2)$. Suppose that $\xi = \langle \boldsymbol{\zeta}, \widehat{\W}_T\rangle$ and $f\equiv 0$. Let $V^{\bone,\btwo\btwo}$ be the collection of all $\bw \in V$ that can be written as concatenations of $\bone$ and $\btwo\btwo$ blocks including the empty word. For $\bw \in V^{\bone,\btwo\btwo}$, let $N(\bw)$ be the number of such blocks, and let $Q(\bw)$ be the number of $\btwo\btwo$-blocks in $\bw$. Let
\[
    \bar{\boldsymbol{\eta}} \coloneqq \sum_{\bv\in V}\sum_{\bw\in V^{\bone,\btwo\btwo}} 2^{-Q(\bw)}\boldsymbol{\zeta}^{\bv\bw} \sum_{k=0}^{N(\bw)} (-1)^k\frac{T^{N(\bw)-k}}{(N(\bw)-k)!} \big(\bone^{\otimes k}\shuprod\bv\big).
\]
Then $(Y,Z)$ given by
\[
    Y_t \coloneqq \langle \bar{\boldsymbol{\eta}}, \widehat{\W}_t\rangle \; \textnormal{and} \; Z_t \coloneqq  \langle \bar{\boldsymbol{\eta}}|_\btwo, \widehat{\W}_t\rangle, \; t \in [0,T],
\]
is the unique square-integrable solution to the martingale representation problem \eqref{eq_bsde_lipschitz}.
\end{proposition}

\begin{remark}
    The tensors $\bar{\boldsymbol{\eta}}$ and $\bar{\boldsymbol{\eta}}|_\btwo$ representing the square-integrable solution to the martingale representation problem \eqref{eq_bsde_lipschitz} are unique in $T(\R^2)$ by \textnormal{\citeauthor*{aqsha2026solving} \cite[Lemma~3.9]{aqsha2026solving}}.
\end{remark}

\begin{proof}
For every $\bw\in V$,
\[
    \cE^\bw(t,T) =
    \begin{cases}
    \displaystyle \frac{2^{-Q(\bw)}}{N(\bw)!}(T-t)^{N(\bw)}, & \bw\in V^{\bone,\btwo\btwo}, \\[1em]
    0, & \bw\notin V^{\bone,\btwo\btwo}.
    \end{cases}
\]
Let $\boldsymbol{\eta}$ be given by \eqref{eq_eta_alpha_beta}. Then
\[
    \langle \boldsymbol{\eta}(t),\widehat{\W}_t\rangle = \sum_{\bv\in V}\sum_{\bw\in V^{\bone,\btwo\btwo}} \boldsymbol{\zeta}^{\bv\bw} \frac{2^{-Q(\bw)}}{N(\bw)!} (T-t)^{N(\bw)} \widehat{\W}^{\bv}_t, \; t \in [0,T].
\]
Using
\[
    \frac{(T-t)^n}{n!} = \sum_{k=0}^n (-1)^k\frac{T^{n-k}}{(n-k)!}\frac{t^k}{k!} \; \textnormal{and} \; \frac{t^k}{k!}=\langle \bone^{\otimes k},\widehat{\W}_t\rangle, \; n \in \N,
\]
and the shuffle identity
\[
    \langle \bone^{\otimes k}\shuprod\bv,\widehat{\W}_t\rangle = \langle \bone^{\otimes k},\widehat{\W}_t\rangle \langle \bv,\widehat{\W}_t\rangle,
\]
we obtain
\[
    \langle \boldsymbol{\eta}(t),\widehat{\W}_t\rangle = \langle \bar{\boldsymbol{\eta}},\widehat{\W}_t\rangle, \; t \in [0,T].
\]
Since
\[
    \langle \boldsymbol{\eta}(t),\widehat{\W}_t\rangle = \langle \boldsymbol{\zeta},\widehat{\W}_T\rangle - \int_t^T \langle \boldsymbol{\eta}(s)|_{\btwo},\widehat{\W}_s\rangle \d W_s, \; t \in [0,T], \; \textnormal{$\P$--a.s.},
\]
and by \Cref{lem_ito_signature}
\[
    \langle \bar{\boldsymbol{\eta}},\widehat{\W}_t\rangle = \langle \boldsymbol{\zeta},\widehat{\W}_T\rangle - \int_t^T \bigg(\Big\langle \bar{\boldsymbol{\eta}}|_{\bone} + \frac{1}{2} \bar{\boldsymbol{\eta}}|_{\btwo\btwo},\widehat{\W}_s\Big\rangle \bigg)\d s - \int_t^T \langle \bar{\boldsymbol{\eta}}|_{\btwo},\widehat{\W}_s\rangle\d W_s, \; t \in [0,T], \; \textnormal{$\P$--a.s.},
\]
we obtain from uniqueness in the decomposition of a special semi-martingale that $\langle \boldsymbol{\eta}(s)|_{\btwo},\widehat{\W}_s\rangle  = \langle \bar{\boldsymbol{\eta}}|_{\btwo},\widehat{\W}_s\rangle$, \textnormal{$\P\otimes\d s$--a.e.,} which yields
\[
    \langle \bar{\boldsymbol{\eta}},\widehat{\W}_t\rangle = \langle \boldsymbol{\zeta},\widehat{\W}_T\rangle - \int_t^T \langle \bar{\boldsymbol{\eta}}|_{\btwo},\widehat{\W}_s\rangle \d W_s, \; t \in [0,T], \; \textnormal{$\P$--a.s.}
\]
This completes the proof.
\end{proof}

\section{An application to pricing and hedging}\label{sec_application}

We give a brief application of the previous results to mathematical finance. For this, we follow the setting in \cite[Sections~1.1 and 2.3]{el1997backward}. Suppose that the prices per unit of a riskless asset $B = (B(t))_{t \in [0,T]}$ and a risky asset $S = (S_t)_{t \in [0,T]}$ satisfy
\[
    \d B(t) =B(t)(-\alpha(t)) \d t \; \textnormal{and} \; \d S_t = S_t \big(b(t) \d t + \sigma(t) \d W_t\big),
\]
with $B(0) \in (0,\infty)$ and, $\P$--a.s., $S_0 \in (0,\infty)$.
Here, the functions $\alpha,b,\sigma : [0,T]\longrightarrow \R$ are bounded and Borel-measurable with $\sigma(t) \neq 0$ for all $t \in [0,T]$, and $\sigma^{-1}$ is bounded.
We denote by $\varphi = (\varphi^0_t,\vartheta_t)_{t \in [0,T]}$ a generic $\F$-progressive investment strategy, where $\varphi^0_t$, resp. $\vartheta_t$, denotes the number of shares held in the riskless asset, resp. in the risky asset, at time $t$. Then $\varphi$ is self-financing if $(\varphi^0,\vartheta)$ is $(B,S)$-integrable, that is,
\[
    \int_0^T \big(|\varphi^0_u\alpha(u)| + |\vartheta_u|^2\big) \d u < \infty, \; \textnormal{$\P$--a.s.},
\]
and
\[
    \varphi^0_t B(t) + \vartheta_t S_t \eqqcolon Y^\varphi_t = Y^\varphi_0 + \int_0^t \varphi^0_u \d B(u) + \int_0^t \vartheta_u \d S_u, \; t \in [0,T], \; \textnormal{$\P$--a.s.}
\]
Writing $Z^\varphi_t \coloneqq \vartheta_t S_t \sigma(t)$ and $\beta(t) \coloneqq -(b(t)+\alpha(t))/\sigma(t)$, we obtain
\begin{equation}\label{eq_value_process_bsde}
    \d Y^\varphi_t = - (\alpha(t) Y^\varphi_t + \beta(t) Z^\varphi_t) \d t + Z^\varphi_t \d W_t,
\end{equation}
which is reminiscent of the dynamics of a BSDE.
Conversely, if an $\F$-progressive pair $(Y,Z)$ is such that $\int_0^T|Z_u|^2\d u < \infty$, $\P$--a.s., and \eqref{eq_value_process_bsde} holds, then $\varphi = (\varphi^0,\vartheta)$ given by
\begin{equation}\label{eq_from_bsde_to_self_financing}
    \varphi^0_t \coloneqq \frac{1}{B(t)}\bigg(Y_t - \frac{Z_t}{\sigma(t)}\bigg) \; \textnormal{and} \; \vartheta_t \coloneqq \frac{Z_t}{\sigma(t)S_t}
\end{equation}
is a self-financing trading strategy with $Y^\varphi_t = Y_t$, $t \in [0,T]$, outside some $\P$--null set.

\medskip
Let $\xi$ denote the payoff of a European option exercisable at time $T$, and let $\pi^\textnormal{s}(\xi)$ be its superhedging price.
By combining \cite[Theorem~1.1 and Remark~b on p.~24]{el1997backward} with \Cref{prop_affine_generator} and \eqref{eq_from_bsde_to_self_financing}, we obtain the following result, which identifies the price and hedging strategy for a payoff that is a polynomial in the signature.

\begin{proposition}
    Suppose $\xi = \langle \boldsymbol{\zeta}, \widehat{\W}_T\rangle$ for $\boldsymbol{\zeta} \in T(\R^2)$, and that $\alpha$ and $\beta$ are bounded and Borel-measurable. Let $\boldsymbol{\zeta}(t)\coloneqq \boldsymbol{\zeta}|_{\cE(t,T)}$, $t \in [0,T]$, where $\cE(t,T)$ is defined as in \textnormal{\Cref{prop_affine_generator}}. Then
    $\varphi = (\varphi^0,\vartheta)$ given by
    \[
        \varphi^0_t = \bigg\langle \frac{1}{B(t)}\bigg(\boldsymbol{\zeta}(t)-\frac{\boldsymbol{\zeta}(t)|_\btwo}{\sigma(t)}\bigg),\widehat{\W}_t \bigg\rangle 
        \; \textnormal{and} \; 
        \vartheta_t = \bigg\langle \frac{\boldsymbol{\zeta}(t)|_\btwo}{\sigma(t)S_t}, \widehat{\W}_t \bigg\rangle, \; t \in [0,T],
    \]
    is a self-financing strategy satisfying, \textnormal{$\P$--a.s.}, $\pi^\textnormal{s}(\xi)=\boldsymbol{\zeta}(0)^\bnothing = Y^\varphi_0$ and $Y^\varphi_T = \xi$.
\end{proposition}

\section{Conclusions}\label{sec_further}

We showed how solutions to a deterministic tensor differential equation yield representations of both components of a BSDE solution driven by Brownian motion. In the case of an affine generator with multiplying deterministic coefficients (relative to $Y$ and $Z$) and an additive signature-dependent term $\langle\boldsymbol{\gamma}(s),\widehat{\W}_s\rangle$, we obtained an explicit solution for terminal conditions that themselves are finite linear combinations of the signature components. We subsequently applied these representations to obtain explicit pricing and hedging formulas.

\medskip
The verification result of \Cref{thm_main_verification} allows further classes of generators to be considered. For example, generators that are affine with signature-dependent coefficients or that are quadratic in $Z$, to name two concrete cases. However, the finite-dimensional construction of solutions to the tensor differential equation above does not extend immediately. This is due to the fact that, in both cases, the tensor differential equation contains shuffle products, whose contributions do not preserve a fixed bound on the tensor degree, in general.

{\small
\bibliography{main_bibliography}

\begin{thebibliography}{19}
\providecommand{\natexlab}[1]{#1}
\providecommand{\url}[1]{\texttt{#1}}
\expandafter\ifx\csname urlstyle\endcsname\relax
  \providecommand{\doi}[1]{doi: #1}\else
  \providecommand{\doi}{doi: \begingroup \urlstyle{rm}\Url}\fi

\bibitem[Abi~Jaber and G{\'e}rard(2025)]{abi2025signature}
E.~Abi~Jaber and L.-A. G{\'e}rard.
\newblock Signature volatility models: pricing and hedging with {F}ourier.
\newblock \emph{SIAM Journal on Financial Mathematics}, 16\penalty0
  (2):\penalty0 606--642, 2025.

\bibitem[Abi~Jaber et~al.(2024)Abi~Jaber, G{\'e}rard, and Huang]{abi2024path}
E.~Abi~Jaber, L.-A. G{\'e}rard, and Y.~Huang.
\newblock Path-dependent processes from signatures.
\newblock \emph{ArXiv preprint arXiv:2407.04956}, 2024.

\bibitem[Abi~Jaber et~al.(2025)Abi~Jaber, Hainaut, and
  Motte]{jaber2025signature}
E.~Abi~Jaber, D.~Hainaut, and E.~Motte.
\newblock Signature approach for pricing and hedging path-dependent options
  with frictions.
\newblock \emph{ArXiv preprint arXiv:2511.23295}, 2025.

\bibitem[Abi~Jaber et~al.(2026{\natexlab{a}})Abi~Jaber, Attal, and
  Sotnikov]{abi2026stochastic}
E.~Abi~Jaber, E.~Attal, and D.~Sotnikov.
\newblock Stochastic control with signatures via {R}iccati equations on the
  tensor algebra.
\newblock \emph{ArXiv preprint arXiv:2607.03986}, 2026{\natexlab{a}}.

\bibitem[Abi~Jaber et~al.(2026{\natexlab{b}})Abi~Jaber, Attal, and
  Sotnikov]{jaber2026affine}
E.~Abi~Jaber, E.~Attal, and D.~Sotnikov.
\newblock Affine structure of the brownian signature.
\newblock \emph{ArXiv preprint arXiv:2606.29622}, 2026{\natexlab{b}}.

\bibitem[Abi~Jaber et~al.(2026{\natexlab{c}})Abi~Jaber, Rey, and
  Sotnikov]{jaber2026malliavin}
E.~Abi~Jaber, C.~Rey, and D.~Sotnikov.
\newblock Malliavin calculus for signatures with applications to finance.
\newblock \emph{ArXiv preprint arXiv:2604.22528}, 2026{\natexlab{c}}.

\bibitem[Aqsha et~al.(2026)Aqsha, Bank, and
  S{\'a}nchez-Betancourt]{aqsha2026solving}
A.~Aqsha, P.~Bank, and L.~S{\'a}nchez-Betancourt.
\newblock Solving linear-quadratic stochastic control problems with signatures.
\newblock \emph{ArXiv preprint arXiv:2602.23473}, 2026.

\bibitem[Bao et~al.(2011)Bao, Delbaen, and Hu]{bao2011backward}
X.~Bao, F.~Delbaen, and Y.~Hu.
\newblock Backward {SDE}s with superquadratic growth.
\newblock \emph{Probability Theory and Related Fields}, 150\penalty0
  (1--2):\penalty0 145--192, 2011.

\bibitem[Bayer et~al.(2026)Bayer, dos Reis, Horvath, and
  Oberhauser]{bayer2026signature}
C.~Bayer, G.~dos Reis, B.~Horvath, and H.~Oberhauser.
\newblock \emph{Signature methods in finance---an introduction with
  computational applications}.
\newblock Springer finance. Springer Cham, 2026.

\bibitem[Chen(1957)]{chen1957integration}
K.-T. Chen.
\newblock Integration of paths, geometric invariants and a generalized
  baker--hausdorff formula.
\newblock \emph{Annals of Mathematics}, 65\penalty0 (1):\penalty0 163--178,
  1957.

\bibitem[Chen(1977)]{chen1977iterated}
K.-T. Chen.
\newblock Iterated path integrals.
\newblock \emph{Bulletin of the American Mathematical Society}, 83\penalty0
  (5):\penalty0 831--879, 1977.

\bibitem[Cuchiero et~al.(2026)Cuchiero, Svaluto-Ferro, and
  Teichmann]{cuchiero2026signature}
C.~Cuchiero, S.~Svaluto-Ferro, and J.~Teichmann.
\newblock Signature {SDE}s from an affine and polynomial perspective.
\newblock \emph{Electronic Journal of Probability}, 31\penalty0 (106):\penalty0
  1--57, 2026.

\bibitem[El~Karoui et~al.(1997)El~Karoui, Peng, and Quenez]{el1997backward}
N.~El~Karoui, S.~Peng, and M.-C. Quenez.
\newblock Backward stochastic differential equations in finance.
\newblock \emph{Mathematical Finance}, 7\penalty0 (1):\penalty0 1--71, 1997.

\bibitem[Kallenberg(2021)]{kallenberg2021foundations}
O.~Kallenberg.
\newblock \emph{Foundations of modern probability}, volume~99 of
  \emph{Probability theory and stochastic modelling}.
\newblock Springer, Cham, 3rd edition, 2021.

\bibitem[Lyons(1998)]{lyons1998differential}
T.~J. Lyons.
\newblock Differential equations driven by rough signals.
\newblock \emph{Revista Matem\'{a}tica Iberoamericana}, 14\penalty0
  (2):\penalty0 215--310, 1998.

\bibitem[Pardoux and Peng(1990)]{pardoux1990adapted}
\'E. Pardoux and S.~Peng.
\newblock Adapted solution of a backward stochastic differential equation.
\newblock \emph{System and Control Letters}, 14\penalty0 (1):\penalty0 55--61,
  1990.

\bibitem[Pardoux and Peng(1992)]{pardoux1992backward}
\'E. Pardoux and S.~Peng.
\newblock Backward stochastic differential equations and quasilinear parabolic
  partial differential equations.
\newblock In B.L. Rozovskii and R.B. Sowers, editors, \emph{Stochastic partial
  differential equations and their applications. Proceedings of IFIP WG 7/1
  international conference University of North Carolina at Charlotte, NC June
  6--8, 1991}, volume 176 of \emph{Lecture notes in control and information
  sciences}, pages 200--217. Springer, 1992.

\bibitem[Xing(2012)]{xing2012on}
H.~Xing.
\newblock On backward stochastic differential equations and strict local
  martingales.
\newblock \emph{Stochastic Processes and their Applications}, 122\penalty0
  (6):\penalty0 2265--2291, 2012.

\bibitem[Zhang(2017)]{zhang2017backward}
J.~Zhang.
\newblock \emph{Backward stochastic differential equations---from linear to
  fully nonlinear theory}, volume~86 of \emph{Probability theory and stochastic
  modelling}.
\newblock Springer-Verlag New York, 2017.

\end{thebibliography}
}

\end{document}